\documentclass[10pt]{amsart}
\usepackage[cp1251]{inputenc}
\usepackage[english,russian]{babel}
\usepackage{amsmath}
\usepackage{amssymb}
\usepackage{amsfonts}

\begin{document}
\sloppy
\begin{center}
{\large\bf $(n+1)$-ARY DERIVATIONS OF SIMPLE MALCEV ALGEBRAS}\\

\hspace*{8mm}

{\large\bf Ivan Kaygorodov}\\
e-mail: kib@math.nsc.ru

{\it 
Sobolev Inst. of Mathematics\\ 
Novosibirsk, Russia\\}
\end{center}

\

\medskip

\

\begin{center} {\bf Abstract: }\end{center}                                                                    
We defined $(n+1)$-ary derivations of $n$-ary algebras. 
We described 3-derivations of simple Malcev algebras
and 4-ary derivations of simple ternary Malcev algebra $M_8$.

\medskip


\medskip

\section{Введение}
Одним из способов обобщения дифференцирований является $\delta$-диф\-фе\-рен\-ци\-рование. 
Под $\delta$-диф\-ференцированием алгебры $A$, при $\delta$ --- фиксированном элементе основного поля,  мы понимаем линейное отображение 
$\phi: A\rightarrow A,$ такое что для произвольных $x,y \in A$ верно
$$\phi(xy)=\delta(\phi(x)y+x\phi(y)).$$
В свое время, $\delta$-дифференцирования изучались в работах 
\cite{Fil}-\cite{kay_nary}, где были описаны $\delta$-дифференцирования 
первичных лиевых \cite{Fil,Fill}, первичных альтернативных и мальцевских \cite{Filll} алгебр, 
простых \cite{kay_lie, kay_lie2} и первичных \cite{Zus} лиевых супералгебр, 
полупростых конечномерных йордановых алгебр \cite{kay,kay_lie2} и супералгебр \cite{kay, kay_lie2,kay_zh,kay_de5}, 
алгебр Филиппова малых размерностей и простых конечномерных алгебр Филиппова \cite{kay_nary}, 
а также простой тернарной алгебры Мальцева $M_8$ \cite{kay_nary}. 
В частности, были построены примеры нетривиальных $\delta$-дифференцирований 
для некоторых алгебр Ли \cite{Fill,Zus, kay_gendelta}, простых йордановых супералгебр \cite{kay_zh, kay_de5} и некоторых $n$-арных алгебр Филиппова \cite{kay_nary}.

В тоже время, $\delta$-дифференцирование является частным случаем квазидифференцирования и обобщенного дифференцирования. 
Под обобщенным дифференцированием $D$ мы понимаем такое линейное отображение, что существуют линейные отображения $E$ и $F,$
связанные с $D$ условием, таким что для произвольных $x,y \in A$ верно
$$D(xy)=E(x)y+xF(y).$$
Если вдобавок к этому $E=F,$ то $D$ --- квазидифференцирование. 
Тройки $(D,E,F),$ где $D$ --- обобщенное дифференцирование, а $E,F$ --- связанные с ним линейные отображения, называются тернарными дифференцированиями. 
Квазидифференцирования, обобщенные дифференцирования и тернарные дифференцирования рассматривались в работах \cite{LL}-\cite{kg_naryfil}.
В частности, изучались обобщенные и тернарные дифференцирования 
алгебр Ли \cite{LL},
супералгебр Ли \cite{ZhangRY}, 
ассоциативных алгебр \cite{komatsu_nak_03, komatsu_nakajima04}, 
обобщенных алгебр Кэли-Диксона \cite{JG3},
йордановых алгебр \cite{shest}
и алгебр Филиппова \cite{kg_naryfil}.

Понятие тернарного дифференцирования для бинарной алгебры допускает обобщение на случай $n$-арных алгебр. 
В данном случае, под $(n+1)$-арным дифференцированием $n$-арной алгебры $A$ мы подразумеваем такой набор $(f_0,f_1, \ldots, f_n) \in End(A)^{n+1},$
что для произвольных $x_1, \ldots, x_n \in A$ верно 
$$f_0[x_1, \ldots, x_n]=\sum\limits_{i=1}^{n}[x_1, \ldots, f_i(x_i), \ldots , x_n].$$
Соответственно, для $(n+1)$-арного квазидифференцирования необходимо дополнительно требовать $f_1=f_2=\ldots=f_n.$ 
Ясно, что если $\psi_1, \ldots, \psi_n$ --- элементы центроида $n$-арной алгебры, 
а $D$ --- дифференцирование $n$-арной алгебры, то 
наборы $(\sum \psi_i, \psi_1, \ldots, \psi_n)$ и $(D, D, \ldots, D)$ --- являются $(n+1)$-арными дифференцированиями. 
Приведенные два вида $(n+1)$-арных дифференцирований, а также их линейные комбинации, мы будем считать тривиальными. 
Ясно, что наибольший интерес представляет вопрос нахождения $(n+1)$-арных дифференцирований, отличных от тривиальных.
Легко заметить следующее 

\medskip

{\bf  Утверждение.} Пусть $A$ --- (анти)коммутативная $n$-арная алгебра и $(f_0, f_1, \ldots, f_n)$ --- $(n+1)$-арное дифференцирование, 
тогда для любой подставновки $\sigma \in S_n$ верно, что $(f_0, f_{\sigma(1)}, \ldots, f_{\sigma(n)})$ --- $(n+1)$-арное дифференцирование. 

\medskip

{\bf Доказательство.} Доказательство данного факта вытекает из известного утверждения о разложении произвольной подстановки в произведение 
транспозиций и очевидного утверждения леммы если $\sigma$ --- транспозиция. Утверждение доказано.

\medskip

Обозначим пространство дифференцирований, $\delta$-дифференцирований, квазидифференцирований и обобщенных дифференцирований, 
соответственно через $Der(A), Der_{\delta}(A), QDer(A)$ и $GDer(A).$
Очевидно, что мы имеем цепочку включений 
\begin{eqnarray}\label{vklyach}
Der(A) \subseteq Der_{\delta}(A) \subseteq QDer(A) \subseteq GDer(A) \subseteq  End(A).
\end{eqnarray}
Стоит отметить, что если $n$-арная алгебра $A$ с ненулевым умножением и характеристика поля отлична от $n-1$, 
то первое включение всегда строгое. 
В противном случае, в силу того, что тождественное отображение явлется элементом центроида и $\frac{1}{n}$-дифференцированием,
мы получили бы противоречение с определением умножения в алгебре $A$. Понятно, что в случае бинарных алгебр ограничение на 
характеристику поля является не существенным.
Ясно, что каждая из алгебр $Der(A), Der_{\delta}(A), QDer(A)$ и $GDer(A)$ относительно коммутаторного умножения становится алгеброй Ли. 
Пусть $Ann(GDer(A))$ --- аннулятор алгебры Ли обобщенных дифференцирований алгебры $A$.
Отметим, что $Ann(GDer(A))$ не тривиален. 
Действительно, там лежат отображения вида $\alpha \cdot id,$ где $\alpha$ --- элемент основного поля.
Обозначим $$\Delta(A)=GDer(A)/Ann(GDer(A)).$$ 
Также нас будет интересовать структура алгебры Ли $\Delta(A).$

Отметим, что в работе \cite{kay_nary} было показано, что простая 8-мерная тернарная алгебра Мальцева $M_8$ не имеет нетривиальных $\delta$-дифференцирований.
Исследованию $(n+1)$-арных дифференцирований была посвящена работа \cite{kg_naryfil}, где были описаны $(n+1)$-арные и обобщенные дифференцирования 
полупростых конечномерных алгебр Филиппова над алгебраически замкнутым полем характеристики нуль. 

Основной целью данной работы, результаты которой аннонсированы в \cite{kay_nplus},
является исследование обобщенных и $(n+1)$-арных дифференцирований 
простых $n$-арых алгебр Мальцева, в случае бинарной алгебры Мальцева $M_7$
(и, следовательно, также любой простой нелиевой бинарной алгебры Мальцева), и тернарной алгебры Мальцева $M_8$.
В работе дано полное описание квазидифференцирований, обобщенных дифференцирований и $(n+1)$-арных дифференцирований алгебр $M_7$ и $M_8$.
В итоге, с использованием результатов \cite{Filll,kay_nary}, мы имеем, что для алгебр $M_7$ и $M_8$
 цепочка включений (\ref{vklyach}) принимает вид 
\begin{eqnarray*}
Der(A) \subset Der_{\delta}(A) = QDer(A) = GDer(A) \subset  End(A).
\end{eqnarray*}

\section{$3$-арные дифференцирования бинарной алгебры Мальцева $M_7$.} 

Алгебры Мальцева, возникающие в \cite{mal_loops}, являются естественным обобщением алгебр Ли и удовлетворяют соотношениям
$$J(x,y,xz)=J(x,y,z)x, \mbox{ где } J(x,y,z)=(xy)z+(yz)x+(zx)y.$$

Хорошо известно \cite{kuz_sh}, 
что существует только одна простая нелиева алгебра Мальцева, 
которой является семимерная алгебра $M_7,$ получающаяся из алгебры Кэли---Диксона.

В случае алгебраически замкнутого поля $P$ характеристики отличной от 2 и 3, 
в алгебре $M_7$ можно выбрать следующий базис --- так называемый расщепляемый базис 
$$B_{M_7}=\{ h,x,y,z,x',y',z'\}$$
с таблицей умножения

$$hx=2x, hy=2y, hz=2z,$$
$$hx'=-2x',hy'=-2y',hz'=-2z',$$
$$xx'=yy'=zz'=h,$$
$$xy=2z',yz=2x',zx=2y',$$
$$x'y'=-2z,y'z'=-2x,z'x'=-2y,$$
где все отсутствующие произведения равны нулю (см. \cite{kuz_sh}).
Каждая тройка элементов $\{h,x,x'\}, \{h,y,y'\},\{h,z,z'\}$ образует стандартный базис трехмерной простой алгебры Ли,
а подпространство $P\cdot h$ является картановской подалгебой алгебры $M_7$.

\medskip

{\bf Теорема 1.} 
Простая алгебра Мальцева $M_7$ не имеет нетривиальных $3$-арных дифференцирований.

\medskip

{\bf Доказательство.} 
Пусть $(D,E,F)$ --- тернарное дифференцирование алгебры $M_7$.
Положим $$D(h)=\gamma h + \sum\limits_{v \in B_{M_7}, v \neq h} \gamma_v v.$$
В этих обозначениях, $(\gamma \cdot id, \frac{\gamma}{2}\cdot id, \frac{\gamma}{2}\cdot id)$ --- тернарное дифференцирование.
В дальнейшем доказательстве теоремы, тернарное дифференцирование 
$$(D-\gamma \cdot id, E-\frac{\gamma}{2} \cdot id, F-\frac{\gamma}{2} \cdot id)$$
мы будем обозначать $(D,E,F).$

Пусть $$W(v)=W^v_hh+W^v_xx+W^v_yy+W^v_zz+W^v_{x'}x'+W^v_{y'}y'+W^v_{z'}z',$$
где $W \in \{D,E,F\}, v \in B_{M_7}, W_w^v \in P.$ 

Заметим, что 
$$D(h)=D(xx')=E(x)x'+xF(x')=$$
$$E_x^xh-2E_h^{x}x'+2E_{y'}^xz-2E_{z'}^xy-2F_h^{x'}x+2F_y^{x'}z'-2F_z^{x'}y'+F_{x'}^{x'}h.$$
Но, в тоже время,
$$D(h)=-D(x'x)=-E(x')x-x'F(x)=F(x)x'+xE(x')=$$
$$F_x^xh-2F_h^{x}x'+2F_{y'}^xz-2F_{z'}^xy-2E_h^{x'}x+2E_y^{x'}z'-2E_z^{x'}y'+E_{x'}^{x'}h.$$
Сравнивая два полученных равенства, мы имеем 
\begin{eqnarray}
\label{mal1-0}
F_{y'}^x=E_{y'}^x,F_{z'}^x=E_{z'}^x,E_x^x+F_{x'}^{x'}=E_{x'}^{x'}+F_x^x=0.
\end{eqnarray}
Отметим, что $0=D(xx)=E(x)x+xF(x),$ то есть $$E_y^x=F_y^x,E_z^x=F_z^x,E_h^x=F_h^x.$$
В тоже время, $0=D(hh)=E(h)h+hF(h),$ откуда $$E^h_w=F^h_w, w\in  B_{M_7} \setminus \{h\}.$$
Тогда легко видеть, что 
$$2D(x)=D(hx)=E(h)x+hF(x)=$$
$$E_h^hx-2E_y^hz'+2E^h_zy'-E^h_{x'}h+2F_y^xy+2F_z^xz-2F_{x'}^xx'-2F_{y'}^xy'-2F_{z'}^xz'+F_x^xx,$$
то есть 
\begin{eqnarray}
 D_y^x=F_y^x, D_z^x=F_z^x. 
\end{eqnarray}

Заметим, что 
$$-2D(x)=D(y'z')=E(y')z'+y'F(z')=$$
\begin{eqnarray}
\label{mal_2.0} 
E_{y'}^{y'}x-2E_h^{y'}z'+E_z^{y'}h+2E_{x'}^{y'}y+2F_h^{z'}y'-F_y^{z'}h+2F_{x'}^{z'}z+F_{z'}^{z'}x,
\end{eqnarray}
то есть 
\begin{eqnarray}
\label{mal_2.1}
2D_h^x=F_y^{z'}-E_z^{y'}, D_y^x=-E_{x'}^{y'}, \\ 
D_z^x=-F_{x'}^{z'}, D_{z'}^x=E_{h}^{y'}, D_{x'}^x=0,D_{y'}^x=-F_{h}^{z'}. 
\label{mal_2}\end{eqnarray}
Также можно заметить, что 
$$-2D_x^x=E_{y'}^{y'}+F_{z'}^{z'}=-F_y^y-E_z^z=2D_{x'}^{x'},$$
где первое равенство вытекакает из (\ref{mal_2.0}), второе из (\ref{mal1-0}), а третье из (\ref{mal_2.0}) путем замены $w$ на $w',$
где $w \in \{x,y,z \},$  в вычислениях.

Ясно, что аналогично мы можем получить $$D_z^z=-D_{z'}^{z'}\mbox{ и }D_y^y=-D_{y'}^{y'}.$$
Легко видеть, что отображение $D^*$, заданное по правилу 
$$D^*(h)=0, D^*(w)=G_w^ww, G_w^w \in P, w \in B_{M_7} \setminus \{h\}, \mbox{ где }$$
$$G_z^z+G_y^y+G_x^x=0, G_w^w=-G_{w'}^{w'}, w \in \{x,y,z\},$$
является дифференцированием алгебры $M_7$. 
Поэтому, в дальнейшем мы можем рассматривать тернарное дифференцирование $(D,E,F)-(D^*,D^*,D^*)$, 
которое будем обозначать как $(D,E,F).$
Таким образом, 
$$D_w^w=0, w \in B_{M_7}\setminus \{x\}, D_x^x \neq 0.$$
Отметим, что 
$$0=D(xy')=E(x)y'+xF(y')=$$
$$-2E_h^xy'+E_y^xh-2E_{x'}^xz+2E_{z'}^xx-2F_h^{y'}x+2F_y^{y'}z'-2F_z^{y'}y'+F_{x'}^{y'}h,$$
что дает 
\begin{eqnarray}\label{mal_3}
E_{x'}^x=0, F_y^{y'}=0, E_h^x=-F_z^{y'}.
\end{eqnarray}
Ясно, что в данных рассуждениях мы могли вместо пары $(x,y')$ брать пары вида $(v,w')$ и $(w',v)$, где $w\neq v$ и $w, v \in \{ x,y,z\}.$
Таким образом, мы получаем
$$F_{w}^{w'}=0, E_{w}^{w'}=0, w \in B_{M_7} \setminus \{h\} \mbox{ и } w''=w,$$
\begin{eqnarray}
\label{mal_4} F_y^{z'} &=&E_h^x, \\
\label{mal_4.1} F_h^{z'}&=&-E_{y'}^x.
\end{eqnarray}

Докажем, что $D_h^x=E_h^x.$ Для этого отметим, что выполняется следующая цепочка равенств 
$$2D_h^x={F_y^{z'}-E_z^{y'}}={E_h^x+F_h^x}=2E_h^x.$$
Первое равенство следует из (\ref{mal_2.1}). Второе равенство вытекает из (\ref{mal_3}) и (\ref{mal_4}).

Покажем, что $D_{y'}^x=E_{y'}^x.$
Для этого достаточно заметить, что 
$$D_{y'}^x=-F_h^{z'} \mbox{ и }F_h^{z'}=-E_{y'}^x.$$
Первое равенство есть последнее равенство в (\ref{mal_2}), а второе --- соотношение (\ref{mal_4.1}).

Таким образом, мы показали, что $D^x_w=E^x_w=F^x_w, w \in B_{M_7} \setminus \{x\}.$
Легко заметить, что верен и аналогичный результат для элементов $B_{M_7} \setminus \{x \}$, 
то есть $D^u_w=E^u_w=F^u_w, w \in B_{M_7} \setminus \{u\}.$
Ясно, что $D$ мы можем представить в виде суммы дифференцирования $D_*$ и линейного отображения $\mu$,
такого что 
$$\mu(x')=-\mu_x x', \mu(x)=\mu_xx, \mu(w)=0, w \in B_{M_7} \setminus \{x, x'\}.$$
Тогда $(\mu,E-D_*,F-D_*)$ является тернарным дифференцированием, которое мы обозначим  
$(\mu,\psi,\chi),$
причем $\psi$ и $\chi$ на элементах базиса алгебры $M_7$ действуют скалярно, то есть 
$$\psi(w)=\psi_ww,\chi(w)=\chi_ww \mbox{ для } w\in B_{M_7}.$$
Покажем, что 
$$(\mu,\psi,\chi)=(0,\sigma \cdot id, -\sigma\cdot id), \sigma \in P.$$
Достаточно заметить, что 
$$\mu(wv)=\psi(w)v+w\chi(v), w,v \in B_{M_7}.$$
Откуда легко извлечь, что 
$$-\psi_h=\chi_y=\chi_z=\chi_{y'}=\chi_{z'},$$
$$-\chi_h=\psi_y=\psi_z=\psi_{y'}=\psi_{z'},$$
$$\psi_z=-\chi_{y'},$$
то есть $$\psi_h=\psi_y=\psi_z=\psi_{y'}=\psi_{z'}=-\chi_h=-\chi_z=-\chi_y=-\chi_{z'}=-\chi_{y'}.$$
Отсюда заключаем, что $\mu_x=\psi_{z'}+\chi_{y'}=0.$
Пользуясь логикой предыдущих рассуждений, получаем, что $$\psi_h=\psi_x=\psi_{x'}=-\chi_x=-\chi_{x'}.$$
Что и влечет требуемое.
Таким образом, любое тернарное дифференцирование алгебры $M_7$ представимо в виде суммы следующих тернарных дифференцирований
$$(D^*,D^*,D^*) \mbox{ и }((\alpha+\beta) \cdot id, \beta \cdot id, \alpha \cdot id),$$
где $D^*$ --- дифференцирование алгебры $M_7$ и $\alpha, \beta \in P.$

Для полного доказательства теоремы осталось установить структуру $\Delta(M_7)$.
Хорошо известно \cite{kuz_sh}, 
что алгебра дифференцирований простой семимерной алгебры Мальцева изоморфна алгебре $B_3.$ 
Таким образом, мы получили $\Delta(M_7) \cong B_3.$ 
Теорема доказана.

\medskip

Заметим, что существует только одна простая нелиева алгебра Мальцева над алгебраически замкнутым полем характеристики нуль (см. \cite{kuz_sh}),
которая изоморфна алгебре $M_7$.
Отсюда и из теоремы 1 имеем

\medskip

{\bf Следствие 2.} 
Простая нелиева алгебра Мальцева над алгебраически замкнутым полем характеристики нуль 
не имеет нетривиальных $3$-арных дифференцирований.

\medskip

Заметим, что полупростая конечномерная алгебра Мальцева $A$ над алгебраически замкнутым полем характеристики нуль представляется в виде прямой суммы
простых идеалов $I_k$.
Тогда, следуя схеме доказательств теорем об описании $\delta$-дифференцирований и $3$-арных дифференцирований полупростых йордановых алгебр
(см. \cite{kay,shest}), 
мы можем получить что каждая компонента $3$-арного дифференцирования инвариантна на прямом слагаемом $I_k.$
В свою очередь, согласно \cite{Fill} существуют простые алгебры Ли обладающие нетривиальными 
$\delta$-дифференцированиями 
(а, соответственно, и нетривиальными $3$-арными дифференцированиями), 
например, трехмерная алгебра $sl_2$ обладает нетривиальным $(-1)$-дифференцированием.
Понятно, что алгебра 
$$\underbrace{sl_2 \oplus \ldots \oplus sl_2}\limits_{m \mbox{ слагаемых}} \oplus  \underbrace{M_7 \oplus \ldots \oplus M_7}\limits_{l \mbox{ слагаемых}}$$ 
является полупростой нелиевой алгеброй Мальцева и обладает нетривиальным $3$-арным дифференцированияем.
В частности, 10-мерная нелиева алгебра Мальцева $sl_2 \oplus M_7.$
Из приведенных рассуждений вытекает

\medskip

{\bf Следствие 3.} 
Над алгебраически замкнутым полем характеристики нуль существуют полупростые нелиевые алгебры Мальцева с нетривиальными $3$-арными дифференцированиями.

\medskip

Используя результаты 
А. Н. Гришкова \cite{grish80} об отсутствии простых конечномерных бинарно-лиевых алгебр над алгебраически замкнутым полем характеристики нуль, 
отличных от алгебр Мальцева и алгебр Ли, а также, теорему 1, мы получаем 

\medskip

{\bf Следствие 4.} 
Если простая конечномерная бинарно-лиева алгебра имеет нетривиальные $3$-арные дифференцирования, то она алгебра Ли.

\medskip

\section{$4$-арные дифференцирования тернарной алгебры Мальцева $M_8$.} 

\medskip

В свое время \cite{fil_nar}, В.Т. Филипповым было предложено некое обобщение алгебр Ли на случай $n$-арной операции. 
В последствии, данный класс алгебр получил название алгебры Филиппова.

Класс $n$-арных алгебр Мальцева был определен в \cite{pozh01}, как некоторый естественный класс $n$-арных алгебр, содержащий 
класс $n$-арных алгебр векторного произведения. 
На самом деле, любая алгебра Филиппова является $n$-арной алгеброй Мальцева.
К настоящему времени единственным известным примером простой $n$-арной алгебры Мальцева, не являющейся алгеброй Филиппова, 
служит простая тернарная алгебра Мальцева $M_8$, возникающая на 8-мерной композиционной алгебре. 
В свое время, дифференцирования тернарной алгебры $M_8$ были описаны в работе \cite{pozh06Der}, 
а в работе \cite{pozh05korn} было построено ее корневое разложение и введена структура $\mathbb{Z}_3$-градуировки.

$n$-Арным якобианом мы называем следующую функцию, определенную на $n$-арной алгебре:

\

$J(x_1, \ldots, x_n; y_2, \ldots, y_n)=$

 $$[[x_1, \ldots, x_n], y_2, \ldots y_n]-\sum\limits_{i=1}^{n}[x_1, \ldots, [x_i,y_2, \ldots, y_n], \ldots,x_n].$$

Из определения следует, что если $A$ --- $n$-арная алгебра Филиппова, то 
$$J(x_1,\ldots, x_n; y_2, \ldots, y_n)=0$$
для всех $x_1,\ldots,x_n,y_2,\ldots, y_n \in A$.

$n$-Арной алгеброй Мальцева $(n\geq3)$ мы называем алгебру $L$ с одной антикоммутативной $n$-арной операцией 
$[x_1, \ldots, x_n]$, удовлетворяющей тождеству
$$-J(zR_x,x_2,\ldots,x_n;y_2,\ldots,y_n)=J(z,x_2,\ldots,x_n;y_2,\ldots,y_n)R_x,$$
где $R_x=R_{x_2,\ldots,x_n}$ --- оператор правого умножения: $zR_x=[z,x_2,\ldots,x_n].$

Далее полагаем, что $P$ --- поле характеристики, отличной от 2,3, и обозначаем через $A$ --- композиционную алгебру
над $F$ с инволюцией $a \rightarrow \overline{a}$ и единицей $1$ (см., например, \cite{zsss}). 
Симметрическую билинейную форму $(x,y)=\frac{1}{2}(x\overline{y}+y\overline{x}),$ определенную на $A$, предполагаем невырожденной и через $n(a)$
обозначаем норму элемента $a\in A.$ Определим на $A$ тернарную операцию умножения $[\cdot, \cdot,\cdot]$ правилом
$$[x,y,z]=x\overline{y}z-(y,z)x+(x,z)y-(x,y)z.$$
Тогда $A$ становится тернарной алгеброй Мальцева \cite{pozh01}, которая обозначается через $M(A),$ а если $dim(A)=8,$ то через $M_8.$

Напомним, что дифференцированием тернарной алгебры $M_8$ называются линейные отображения $D$ при произвольных элементах $x,y,z \in M_8$ удовлетворяющие равенству 
$$D[x,y,z]=[D(x),y,z]+[x,D(y),z]+[x,y,D(z)].$$
В работе \cite{pozh06Der} было описаны дифференцирования тернарной алгебры Мальцева $M_8$, 
где было показано, что каждое дифференцирование является внутренним, 
то есть 
$$Der(M_8)=\langle [R_{x,y},R_{x,z}]+R_{x,[y,x,z]} | x,y,z \in M_8 \rangle.$$
Отметим, что в той же работе был построен базис $Der(M_8)$:
\begin{eqnarray*}\label{basis}
\mathbb{B}=\{ 
&&\Delta_{23}-\Delta_{14}, \Delta_{24}+\Delta_{13}, \Delta_{25}-\Delta_{16}, \Delta_{26}+\Delta_{15}, \Delta_{27}+\Delta_{18}, \\
&&\Delta_{28}-\Delta_{17}, \Delta_{34}-\Delta_{12}, \Delta_{35}-\Delta_{17}, \Delta_{36}-\Delta_{18}, \Delta_{37}+\Delta_{15}, \\ 
&&\Delta_{38}+\Delta_{16}, \Delta_{45}-\Delta_{18}, \Delta_{46}+\Delta_{17}, \Delta_{47}-\Delta_{16}, \Delta_{48}+\Delta_{15}, \\ 
&&\Delta_{56}-\Delta_{12}, \Delta_{57}-\Delta_{13}, \Delta_{58}-\Delta_{14}, \Delta_{67}+\Delta_{14}, \Delta_{68}-\Delta_{13}, \Delta_{78}+\Delta_{12} \},
\end{eqnarray*}
где $\Delta_{ij}=e_{ij}-e_{ji}$ и $e_{ij}$ --- обычные матричные единички.


Под $4$-арным дифференцированием тернарной алгебры $M_8$ мы подразумеваем четверку $(D,E,F,G) \in End(A)^4$, такую что для 
произвольных элементов $x,y,z \in M_8$ верно
$$D[x,y,z]=[E(x),y,z]+[x,F(y),z]+[x,y,G(z)].$$

Пусть $1,a,b,c$ --- ортонормированные вектора из $A$. Выберем следующий базис в $A$:
$$\wp=\{e_1=1,e_2=a,e_3=b,e_4=ab,e_5=c,e_6=ac,e_7=bc,e_8=abc\}.$$

Известно, что для каждого $i \in \{ 2, \ldots ,8 \}$ возможно выбрать $j,k,l,m,s,t$, зависящие от $i$, такие, что 
\begin{eqnarray}\label{m8_2}
e_i=e_je_k=e_le_m=e_se_t,\\
e_j=e_se_m=e_ke_i=e_te_l,\\
e_k=e_ie_j=e_me_t=e_se_l,\\
e_l=e_me_i=e_ke_s=e_je_t,\\
e_m=e_ie_l=e_te_k=e_je_s,\\
e_s=e_le_k=e_te_i=e_me_j,\\
\label{m8_3}e_t=e_ie_s=e_ke_m=e_le_j.
\end{eqnarray}

\medskip 

\medskip 

\textbf{Теорема 5.} 
Простая тернарная алгебра Мальцева $M_8$ не имеет нетривиальных $4$-арных дифференцирований.

\

{\bf Доказательство.} 
Пусть $(D,E,F,G)$ --- 4-арное дифференцирование алгебры $M_8$.
Оператор правого умножения $R_{x,y}$ называется регулярным, если в фиттинговом разложении $M=M_0 \oplus M_1$ 
относительно $R_{x,y}$ размерность $M_0$ минимальна \cite{pozh05korn}. 
Согласно \cite[Теорема 1]{pozh05korn}, мы имеем корневое разложение алгебры $M_8:$
$M=M_0 \oplus M_{\alpha} \oplus M_{-\alpha},$ где $\alpha \in F$ такой, что $vR_{x,y}=\pm\alpha v$ для любого $v\in M_{\pm\alpha}.$ 
Также из \cite[Лемма 3]{pozh05korn} известно, что на тернарной алгебре $M_8$ существует нетривиальная градуировка.
Если мы обозначим $M_{\pm\alpha}$ через $M_{\pm1},$ то $$[M_i, M_j, M_k] \subseteq M_{i+j+k(mod3)}.$$

Заметим, что согласно \cite[Лемма 1]{pozh05korn}, 
операторы $R_{e_p,e_q}$ при $p \neq q$ являются регулярными, следовательно, верно
$$D[e_p,e_p,e_q]=[E(e_p),e_p,e_q]+[e_p,F(e_p),e_q]+[e_p,e_p,G(e_q)],$$
то есть $$[E(e_p),e_p,e_q]-[F(e_p),e_p,e_q]=0.$$
Последнее, в силу регулярности оператора $R_{e_p,e_q}$ и $\mathbb{Z}_3$-градуировки $M_8$,
 влечет, что $E(e_p)-F(e_p)= \xi_p e_p$. Аналогичное верно и для пары отображений $F$ и $G$.
Таким образом, мы можем считать, что 
\begin{eqnarray}\label{feg}
E_p^r=F_p^r=G_p^r,\mbox{ при } p \neq r.
\end{eqnarray}

Благодаря тому, что 
$$D(e_s)=D[e_i,e_j,e_l]=[E(e_i),e_j,e_l]+[e_i,F(e_j),e_l]+[e_i,e_j,G(e_l)],$$
выполнив соответствующие операции умножений, мы имеем
\begin{eqnarray}
\label{metka}\sum\limits_{p=1}^{8} D_s^p e_p= 
E_i^1e_t+E_i^ie_s+E_i^ke_m-E_i^me_k-E_i^se_i-E_i^te_1+
\end{eqnarray}
\begin{eqnarray*}
&&F_j^1e_m+F_j^je_s-F_j^ke_t-F_j^me_1-F_j^se_j+F_j^te_k-\\
&&G_l^1e_k+G_l^ke_1+G_l^le_s-G_l^me_t-G_l^se_l+G_l^te_m
\end{eqnarray*}
(к примеру, $[e_i,e_l,e_k]=(e_i \overline{e_l})e_k=(e_le_i)e_k=-e_me_k=e_ke_m=e_t$).
Следовательно, легко вытекает $D_s^i=-E_i^s, D_s^j=-F_j^s, D_s^l=-G_l^s.$
Таким образом, пользуясь соотношением (\ref{feg}), мы можем считать, что $D_s^p=-E_p^s,$ где $p \in \{i,j,l\}$. 
Произвольность индекса $i$ и соотношения (\ref{m8_2}-\ref{m8_3})
позволяют нам сделать вывод, что 
\begin{eqnarray}
\label{defg} D_q^p=-E_p^q=-F_p^q=-G_p^q,\text{ где }p,q \in \{1,\ldots, 8\} \text{ и } p\neq q .
\end{eqnarray}

Пользуясь соотношениями (\ref{defg}) и учитывая (\ref{metka}), мы можем получить
\begin{eqnarray}
\label{DDDD1}D_s^t&=&-D_1^i+D_k^j+D_m^l, \\
D_s^1&=& D_t^i+D_m^j-D_k^l, \\
D_s^m&=&-D_k^i-D_1^j-D_t^l, \\
\label{DDDD2}D_s^k&=& D_m^i-D_t^j+D_1^l.
\end{eqnarray}
Учитывая соотношения (\ref{DDDD1}-\ref{DDDD2}), базис алгебры $Der(M_8)$ и тот факт, что 
$(D^*,D^*,D^*,D^*)$ --- $4$-арное дифференцирование для $D^* \in Der(M_8),$
мы можем считать, что матрица $[D]_{\wp}$ является верхнетреугольной.
Это действительно верно, ведь алгебра $4$-арных дифференцирований является замкнутой относительно сложения 
и, следовательно, мы можем обнулить все коэффициенты в нижнем треугольнике части матрицы $[D]_{\wp}$, 
посредством вычитания соответствующих дифференцирований из отображения $D$.

Полученное $4$-арное дифференцирование мы будем обозначать как и прежде $(D,E,F,G)$ и видим, что 
$$D(e_{\gamma})=\sum\limits_{\beta \leq \gamma} D_{\gamma}^{\beta} e_{\beta}\mbox{ и }
T(e_{\gamma})=\sum\limits_{\gamma \leq \beta} T_{\gamma}^{\beta} e_{\beta},
\mbox{ где }D_{\gamma}^{\beta},T_{\gamma}^{\beta} \in P  \mbox{ и }  T\in \{ E,F,G\}.$$
Легко заметить, что 
$$D[e_6,e_7,e_8]=[E(e_6),e_7,e_8]+[e_6,F(e_7),e_8]+[e_6,e_7,G(e_8)]=(E_6^6+F_7^7+G_8^8)e_5,$$
таким образом, $D(e_5)=D_5^5e_5.$ Аналогично получаем 
$$D[e_5,e_7,e_8]=[E(e_5),e_7,e_8]+[e_5,F(e_7),e_8]+[e_5,e_7,G(e_8)]=(E_5^5+F_7^7+G_8^8)e_6,$$
то есть, $D(e_6)=D_6^6e_6.$ Теперь легко видеть, что 
$$D[e_5,e_6,e_8]=[E(e_5),e_6,e_8]+[e_5,F(e_6),e_8]+[e_5,e_6,G(e_8)]=(E_5^5+F_6^6+G_8^8)e_7,$$
что влечет $D(e_7)=D_7^7e_7.$ Пользуясь предыдущими рассуждениями имеем
$$D[e_6,e_5,e_7]=[E(e_6),e_5,e_7]+[e_6,F(e_5),e_7]+[e_6,e_5,G(e_7)]=(E_6^6+F_5^5+G_7^7)e_8,$$
откуда следует, что $D(e_8)=D_8^8e_8.$ Заметим, что 
$$D_3^1e_1+D_3^2e_2+D_3^3e_3=D(e_3)=D[e_4,e_5,e_6]=$$
$$[E(e_4),e_5,e_6]+[e_4,F(e_5),e_6]+[e_4,e_5,G(e_6)]=(E_4^4+F_5^5+G_6^6)e_3+E_4^7e_8+E_4^8e_7$$
и значит $D(e_3)=D_3^3e_3.$ Осталось заметить, что $D(e_4)=D_4^4e_4$ и $D(e_2)=D_2^2e_2.$
Это вытекает из соотношений $e_4=[e_5,e_3,e_6]$ и $e_2=[e_7,e_4,e_5],$ соответственно.
Таким образом, мы показали, что $D(e_{\gamma})=D_{\gamma} e_{\gamma}$ и, следовательно, 
$$E(e_{\gamma})=E_{\gamma} e_{\gamma},F(e_{\gamma})=F_{\gamma} e_{\gamma},G(e_{\gamma})=G_{\gamma} e_{\gamma}.$$

Пусть $x,y,z$ некоторые различные элементы из базиса $\wp$. Тогда можно заметить, что 
$$[E(x),y,z]+[x,F(y),z]=D[x,y,z]-[x,y,G(z)]=[F(x),y,z]+[x,E(y),z],$$
то есть, можно получить, что $E_x-F_x=E_y-F_y.$ Откуда можно заключить, что 
$E=F+ \alpha \cdot id$ при некотором $\alpha \in M_8.$ Аналогично можно получить, что 
$G=F+ \beta \cdot id$ при некотором $\beta \in M_8.$ Отсюда следует, что 
$$(D,E,F,G)=(D, F+ \alpha \cdot id, F,F+ \beta \cdot id).$$
В силу того, что 
$((\alpha +\beta)\cdot id, \alpha \cdot id, 0,\beta \cdot id)$ --- $4$-арное дифференцирование, 
можем считать, что 
$$(D,E,F,G)=(D, F, F,F).$$
Осталось показать, что $D=3F=3f \cdot id, f\in P.$ 
Для этого заметим, что в силу того, что $(D, F, F,F)$ --- $4$-арное дифференцирование, верно
\begin{eqnarray}
\label{3f1} D_s&=&F_i+F_j+F_l,\\
\label{3f2} D_t&=&F_i+F_k+F_l,\\
\label{3f3} D_t&=&F_i+F_j+F_m,\\
\label{3f4} D_s&=&F_i+F_k+F_m.
\end{eqnarray}
Вычитая из (\ref{3f1}) равенство (\ref{3f2}) и из (\ref{3f3}) равенство (\ref{3f4}), мы получаем
$$D_s-D_t=F_j-F_k=D_t-D_s.$$
Откуда, видим, что $D_t=D_s$ и $F_j=F_k$. Действую аналогично, мы получим $D=3F=3f\cdot id.$ 
Таким образом, мы показали, что произвольное $4$-арное дифференцирование тернарной алгебры Мальцева $M_8$
представимо в виде суммы $4$-арных дифференцирований следующих видов 
$$(D^*,D^*,D^*,D^*) \mbox{ и }(\sum\limits_{i=1}^3\alpha_i\cdot id, \alpha_1 \cdot id, \alpha_2 \cdot id, \alpha_3 \cdot id),$$
где $D^* \in Der(M_8)$ и $\alpha_i \in P$.

Для полного доказательства теоремы нам необходимо отметить, что в \cite{pozh06Der} показана 
структура алгебры дифференцирований тернарной алгебры Мальцева $M_8$.
А именно, доказан изоморфизм с алгеброй $B_3.$ 
Таким образом, легко видеть, что $\Delta(M_8) \cong B_3$.
Теорема доказана.

\medskip 

В результате, из теорем 1 и 5 вытекает

\medskip

\textbf{Следствие 6.}  $\Delta(M_7) \cong \Delta(M_8)$.

\medskip 

В заключение, автор выражает благодарность проф. В. Н. Желябину и проф. А. П. Пожидаеву за внимание к работе и конструктивные замечания.


\end{document}